\documentclass[twoside,12pt]{article}
\usepackage{amsfonts}
\usepackage{CJK}
\usepackage{indentfirst}
\usepackage{bm}
\usepackage{amsmath}

\begin{document}
\begin{CJK*}{GBK}{song}

\begin{center}
\LARGE\bf The Cauchy problem for higher-order linear partial
differential equation
\end{center}

\footnotetext{\hspace*{-.45cm}\footnotesize $^\dag$Corresponding
author. E-mail: yananbiguangqing@sohu.com}

\begin{center}
\rm Guangqing Bi $^{\rm a)\dagger}$, \ \ Yuekai Bi $^{\rm b)}$
\end{center}

\begin{center}
\begin{footnotesize} \sl
${}^{\rm b)}$ School of Electronic and Information Engineering,
BUAA, Beijing 100191, China\\
E-mail: yuekaifly@163.com
\end{footnotesize}
\end{center}

\vspace*{2mm}

\begin{center}
\begin{minipage}{15.5cm}
\parindent 20pt\footnotesize
For the linear partial differential equation
$P(\partial_x,\partial_t)u=f(x,t)$, where
$x\in\mathbb{R}^n,\;t\in\mathbb{R}^1$, with
$P(\partial_x,\partial_t)$ is
$\prod^m_{i=1}(\frac{\partial}{\partial{t}}-a_iP(\partial_x))$ or
$\prod^m_{i=1}(\frac{\partial^2}{\partial{t^2}}-a_i^2P(\partial_x))$,
the authors give the analytic solution of the cauchy problem using
the abstract operators $e^{tP(\partial_x)}$ and
$\frac{\sinh(tP(\partial_x)^{1/2})}{P(\partial_x)^{1/2}}$. By
representing the operators with integrals, explicit solutions are
obtained with an integral form of a given function.
\end{minipage}
\end{center}

\begin{center}
\begin{minipage}{15.5cm}
\begin{minipage}[t]{2.3cm}{\bf Keywords:}\end{minipage}
\begin{minipage}[t]{13.1cm}
cauchy problem, partial differential equation, abstract operators
\end{minipage}\par\vglue8pt
{\bf MSC(2000): } 35G10; 35G05
\end{minipage}
\end{center}

\section{Introduction and main results}
In 1997, Guangqing Bi first introduced the concept of abstract
operators in reference \cite{bi97}, and determined the algorithms of
abstract operators $e^{tP(\partial_x)}$ and
$\frac{\sinh(tP(\partial_x)^{1/2})}{P(\partial_x)^{1/2}}$. Using
this type of operators, in reference \cite{bi99} the author has
obtained the following results:

\textbf{Theorem BI1.} Let $a_1,a_2,\ldots,a_m$ be arbitrary real or
complex numbers different from each other, $P(\partial_x)$ be a
partial differential operator of any order, then $\forall
f(x,t)\in{C^\infty(\mathbb{R}_x^n\times\mathbb{R}_t^1)}$, we have
\begin{equation}\label{I1}
    \left\{\begin{array}{l@{\qquad}l}\displaystyle
    \prod^m_{i=1}(\frac{\partial}{\partial{t}}-a_iP(\partial_x))u=f(x,t),&x\in\mathbb{R}^n,\quad{t}\in\mathbb{R}^1.\\\displaystyle
    \left.\frac{\partial^ju}{\partial{t^j}}\right|_{t=0}=0,&j=0,1,2,\ldots,m-1.
    \end{array}\right.
\end{equation}
\begin{equation}\label{I1'}
u(x,t)=\int^t_0\int^{t-\tau}_0\frac{(t-\tau-\tau')^{m-2}}{(m-2)!}\sum^m_{j=1}\frac{a_j^{m-1}}{\prod^m_{i=1\atop
i\neq{j}}(a_j-a_i)} e^{\tau'a_jP(\partial_x)}f(x,\tau)\,d\tau'd\tau.
\end{equation}

\textbf{Theorem BI2.} Let $a_1,a_2,\ldots,a_m$ be arbitrary real or
complex numbers different from each other, $P(\partial_x)$ be a
partial differential operator of any order, then $\forall
f(x,t)\in{C^\infty(\mathbb{R}_x^n\times\mathbb{R}_t^1)}$, we have
\begin{equation}\label{I2}
    \left\{\begin{array}{l@{\qquad}l}\displaystyle
    \prod^m_{i=1}(\frac{\partial^2}{\partial{t^2}}-a_i^2P(\partial_x))u=f(x,t),&x\in\mathbb{R}^n,\quad{t}\in\mathbb{R}^1.\\\displaystyle
    \left.\frac{\partial^ju}{\partial{t^j}}\right|_{t=0}=0,&j=0,1,2,\ldots,2m-1.
    \end{array}\right.
\end{equation}
\begin{equation}\label{I2'}
u(x,t)=\int^t_0\int^{t-\tau}_0\frac{(t-\tau-\tau')^{2m-3}}{(2m-3)!}\sum^m_{j=1}\frac{a_j^{2m-2}}{\prod^m_{i=1
\atop i\neq{j}}(a_j^2-a_i^2)}
\frac{\sinh(\tau'a_jP(\partial_x)^{1/2})}{a_jP(\partial_x)^{1/2}}f(x,\tau)\,d\tau'd\tau.
\end{equation}

In reference \cite{bi01} Guangqing Bi has obtained the following
results:

\textbf{Theorem BI3.} Let $m\in\mathbb{N}$ and $m>1$,
$P(\partial_x)$ be a partial differential equation of any order.
Then $\forall
f(x,t)\in{C^\infty(\mathbb{R}_x^n\times\mathbb{R}_t^1)}$, we have
\begin{equation}\label{I3}
    \left\{\begin{array}{l@{\qquad}l}\displaystyle
    \left(\frac{\partial^2}{\partial{t^2}}-P(\partial_x)\right)^mu=f(x,t),&x\in\mathbb{R}^n,\quad{t}\in\mathbb{R}^1.\\\displaystyle
    \left.\frac{\partial^ju}{\partial{t^j}}\right|_{t=0}=0,&j=0,1,2,\ldots,2m-1.
    \end{array}\right.
\end{equation}
\begin{equation}\label{I3'}
u(x,t)=\int^t_0\int^{t-\tau}_0\frac{((t-\tau)^2-\tau'^2)^{m-2}}{(2m-2)!!(2m-4)!!}\,
\frac{\sinh\left(\tau'P(\partial_x)^{1/2}\right)}{P(\partial_x)^{1/2}}f(x,\tau)\,d\tau'd\tau.
\end{equation}

By combining the abstract operators and Laplace transform, the
authors have obtained the following results in reference
\cite{bi10}:

\textbf{Theorem BI4.} Let $m\in\mathbb{N}$ and $m>1$,
$P(\partial_x)$ be a partial differential equation of any order.
Then
$\forall{f(x,t)}\in{C^\infty}(\mathbb{R}_x^n\times\mathbb{R}_t^1),\;\varphi_j(x)\in{C^\infty}(\mathbb{R}^n)$,
we have
\begin{equation}\label{I4}
    \left\{\begin{array}{l@{\qquad}l}\displaystyle
    \left(\frac{\partial^2}{\partial{t^2}}-P(\partial_x)\right)^mu=f(x,t),&x\in\mathbb{R}^n,\quad{t}\in\mathbb{R}^1.\\\displaystyle
    \left.\frac{\partial^ju}{\partial{t^j}}\right|_{t=0}=\varphi_j(x),&j=0,1,2,\ldots,2m-1.
    \end{array}\right.
\end{equation}
\begin{eqnarray}\label{I4'}
  u(x,t) &=&
  \int^t_0\int^{t-\tau}_0\frac{((t-\tau)^2-\tau'^2)^{m-2}}{(2m-2)!!\,(2m-4)!!}\,
  \frac{\sinh\left(\tau'P(\partial_x)^{1/2}\right)}{P(\partial_x)^{1/2}}\,f(x,\tau)\,\tau'd\tau'\,d\tau \nonumber\\
   & & +\,\sum^{m-1}_{k=0}(-1)^k{m\choose{k}}P(\partial_x)^k\sum^{2m-1-2k}_{j=0}\frac{\partial^{2m-1-2k-j}}{\partial t^{2m-1-2k-j}}\int^t_0
   \frac{(t^2-\tau^2)^{m-2}\tau}{(2m-2)!!\,(2m-4)!!}\nonumber\\
   & & \times\,\frac{\sinh\left(\tau P(\partial_x)^{1/2}\right)}{P(\partial_x)^{1/2}}\,\varphi_j(x)\,d\tau.
\end{eqnarray}

Using the same method, we have obtained the following results in
this paper:

\textbf{Theorem 1.} Let $a_1,a_2,\ldots,a_m$ be arbitrary real or
complex roots different from each other for
$b_0+b_1x+b_2x^2+\cdots+b_mx^m=0$, and $P(\partial_x,\partial_t)$ be
a partial differential operator defined by
$$P(\partial_x,\partial_t)=\sum^m_{k=0}b_kP(\partial_x)^{m-k}\frac{\partial^k}{\partial{t^k}},
\quad{x}\in\mathbb{R}^n,\;t\in\mathbb{R}^1,\;1<m\in\mathbb{N}.$$
Where $P(\partial_x)$ is a partial differential operator of any
order. Then
$\forall{f(x,t)}\in{C^\infty}(\mathbb{R}_x^n\times\mathbb{R}_t^1),\;\varphi_r(x)\in{C^\infty}(\mathbb{R}^n)$,
we have
\begin{equation}\label{1}
    \left\{\begin{array}{l@{\qquad}l}\displaystyle
    P(\partial_x,\partial_t)u=f(x,t),&x\in\mathbb{R}^n,\quad{t}\in\mathbb{R}^1.\\\displaystyle
    \left.\frac{\partial^ru}{\partial{t^r}}\right|_{t=0}=\varphi_r(x),&r=0,1,2,\ldots,m-1.
    \end{array}\right.
\end{equation}
\begin{eqnarray}\label{1'}
  u(x,t) &=&
  \int^t_0\int^{t-\tau}_0\frac{(t-\tau-\tau')^{m-2}}{(m-2)!}\sum^m_{j=1}\frac{a_j^{m-1}}{\prod^m_{i=1 \atop i\neq{j}}(a_j-a_i)}\,
e^{\tau'a_jP(\partial_x)}f(x,\tau)\,d\tau'd\tau\nonumber\\
   & & +\,\sum^m_{k=1}b_kP(\partial_x)^{m-k}\sum^{k-1}_{r=0}\frac{\partial^{k-1-r}}{\partial t^{k-1-r}}\int^t_0
   \frac{(t-\tau)^{m-2}}{(m-2)!}\nonumber\\
   & & \times\,\sum^m_{j=1}\frac{a_j^{m-1}}{\prod^m_{i=1 \atop i\neq{j}}(a_j-a_i)}\,e^{\tau{a_j}P(\partial_x)}\,\varphi_r(x)\,d\tau.
\end{eqnarray}

\textbf{Theorem 2.} Let $a_1,a_2,\ldots,a_m$ be arbitrary real or
complex roots different from each other, satisfy
$\sum^m_{k=0}b_{2k}x^{2k}=\prod^m_{i=1}(x^2-a_i^2)$, and
$P(\partial_x,\partial_t)$ be a partial differential operators
defined by
$$P(\partial_x,\partial_t)=\sum^m_{k=0}b_{2k}P(\partial_x)^{m-k}\frac{\partial^{2k}}{\partial{t^{2k}}},
\quad{x}\in\mathbb{R}^n,\;t\in\mathbb{R}^1,\;1<m\in\mathbb{N}.$$
Where $P(\partial_x)$ be a partial differential operator of any
order, then
$\forall{f(x,t)}\in{C^\infty}(\mathbb{R}_x^n\times\mathbb{R}_t^1),\;\varphi_r(x)\in{C^\infty}(\mathbb{R}^n)$,
we have
\begin{equation}\label{2}
    \left\{\begin{array}{l@{\qquad}l}\displaystyle
    P(\partial_x,\partial_t)u=f(x,t),&x\in\mathbb{R}^n,\quad{t}\in\mathbb{R}^1.\\\displaystyle
    \left.\frac{\partial^ru}{\partial{t^r}}\right|_{t=0}=\varphi_r(x),&r=0,1,2,\ldots,2m-1.
    \end{array}\right.
\end{equation}
\begin{eqnarray}\label{2'}
  u(x,t) &=&
  \int^t_0\int^{t-\tau}_0\frac{(t-\tau-\tau')^{2m-3}}{(2m-3)!}\sum^m_{j=1}\frac{a_j^{2m-2}}{\prod^m_{i=1 \atop i\neq{j}}(a_j^2-a_i^2)}
\frac{\sinh(\tau'a_jP(\partial_x)^{1/2})}{a_jP(\partial_x)^{1/2}}f(x,\tau)\,d\tau'd\tau\nonumber\\
   & & +\,\sum^m_{k=1}b_{2k}P(\partial_x)^{m-k}\sum^{2k-1}_{r=0}\frac{\partial^{2k-1-r}}{\partial t^{2k-1-r}}\int^t_0
   \frac{(t-\tau)^{2m-3}}{(2m-3)!}\nonumber\\
   & & \times\,\sum^m_{j=1}\frac{a_j^{2m-2}}{\prod^m_{i=1 \atop i\neq{j}}(a_j^2-a_i^2)}
   \frac{\sinh(\tau{a_j}P(\partial_x)^{1/2})}{a_jP(\partial_x)^{1/2}}\,\varphi_r(x)\,d\tau.
\end{eqnarray}

\section{Proof of theorems}

According to the Theorem BI1 and Theorem BI2, we just need to prove
the following corollary of Theorem 1 and Theorem 2:

\textbf{Corollary 1.} Let $a_1,a_2,\ldots,a_m$ be arbitrary real or
complex roots different from each other for
$b_0+b_1x+b_2x^2+\cdots+b_mx^m=0$, and $P(\partial_x,\partial_t)$ be
a partial differential operators defined by
$$P(\partial_x,\partial_t)=\sum^m_{k=0}b_kP(\partial_x)^{m-k}\frac{\partial^k}{\partial{t^k}},
\quad{x}\in\mathbb{R}^n,\;t\in\mathbb{R}^1,\;1<m\in\mathbb{N}.$$
Where $P(\partial_x)$ be a partial differential operator of any
order, then for $\forall\varphi_r(x)\in{C^\infty(\mathbb{R}^n)}$, we
have
\begin{equation}\label{3}
    \left\{\begin{array}{l@{\qquad}l}\displaystyle
    P(\partial_x,\partial_t)u=0,&x\in\mathbb{R}^n,\quad{t}\in\mathbb{R}^1.\\\displaystyle
    \left.\frac{\partial^ru}{\partial{t^r}}\right|_{t=0}=\varphi_r(x),&r=0,1,2,\ldots,m-1.
    \end{array}\right.
\end{equation}
\begin{eqnarray}\label{3'}
  u(x,t)=\sum^m_{k=1}b_kP(\partial_x)^{m-k}\sum^{k-1}_{r=0}\frac{\partial^{k-1-r}}{\partial
t^{k-1-r}}\int^t_0\frac{(t-\tau)^{m-2}}{(m-2)!}
\sum^m_{j=1}\frac{a_j^{m-1}e^{\tau{a_j}P(\partial_x)}}{\prod^m_{i=1
\atop i\neq{j}}(a_j-a_i)}\,\varphi_r(x)\,d\tau.
\end{eqnarray}

\textbf{Proof.} Considering initial conditions, the Laplace
transform of the Eq (\ref{3}) with respect to $t$ is
$$\sum^m_{k=0}b_kP(\partial_x)^{m-k}(s^kU(x,s)-\sum^{k-1}_{r=0}s^{k-1-r}\varphi_r(x))=0,$$
$$\prod^m_{i=1}(s-a_iP(\partial_x))U(x,s)-\sum^m_{k=1}b_kP(\partial_x)^{m-k}\sum^{k-1}_{r=0}s^{k-1-r}\varphi_r(x)=0.$$
Where $U(x,s)=Lu(x,t)$. Let
$G_m(\partial_x,t)=L^{-1}[1/\prod^m_{i=1}(s-a_iP(\partial_x))]$, by
solving $U(x,s)$, we have its inverse Laplace transform:
\begin{eqnarray}\label{3''}
u(x,t)&=& L^{-1}U(x,s)= \sum^m_{k=1}b_kP(\partial_x)^{m-k}\sum^{k-1}_{r=0}L^{-1}\frac{s^{k-1-r}}{\prod^m_{i=1}(s-a_iP(\partial_x))}\,\varphi_r(x)\nonumber\\
   &=& \sum^m_{k=1}b_kP(\partial_x)^{m-k}\sum^{k-1}_{r=0}\frac{\partial^{k-1-r}}{\partial{t^{k-1-r}}}\,G_m(\partial_x,t)\varphi_r(x).
\end{eqnarray}

Now let us solve $G_m(\partial_x,t)$. Considering initial
conditions, the Laplace transform of the Eq (\ref{I1}) with respect
to $t$ is
$$\prod^m_{i=1}(s-a_iP(\partial_x))U(x,s)=F(x,s),\quad{F(x,s)}=Lf(x,t).$$
By solving $U(x,s)$ and using the convolution theorem, we have its
inverse Laplace transform:
$$u(x,t)=L^{-1}U(x,s)=L^{-1}\frac{1}{\prod^m_{i=1}(s-a_iP(\partial_x))}F(x,s)=G_m(\partial_x,t)*f(x,t).$$
By comparing (\ref{I1'}) with $u(x,t)=G_m(\partial_x,t)*f(x,t)$, we
have the expression of the abstract operator $G_m(\partial_x,t)$:
\begin{equation}\label{4}
G_m(\partial_x,t)=\int^t_0\frac{(t-\tau)^{m-2}}{(m-2)!}
\sum^m_{j=1}\frac{a_j^{m-1}}{\prod^m_{i=1 \atop
i\neq{j}}(a_j-a_i)}e^{\tau{a_j}P(\partial_x)}d\tau.
\end{equation}
Applying (\ref{4}) to (\ref{3''}), thus the Corollary 1 is proved.

\textbf{Corollary 2.} Let $a_1,a_2,\ldots,a_m$ be arbitrary real or
complex roots different from each other, which satisfy
$\sum^m_{k=0}b_{2k}x^{2k}=\prod^m_{i=1}(x^2-a_i^2)$, and
$P(\partial_x,\partial_t)$ be a partial differential operator
defined by
$$P(\partial_x,\partial_t)=\sum^m_{k=0}b_{2k}P(\partial_x)^{m-k}\frac{\partial^{2k}}{\partial{t^{2k}}},
\quad{x}\in\mathbb{R}^n,\;t\in\mathbb{R}^1,\;1<m\in\mathbb{N}.$$
Where $P(\partial_x)$ is a partial differential operator of any
order, then $\forall\varphi_r(x)\in{C^\infty(\mathbb{R}^n)}$, we
have
\begin{equation}\label{5}
    \left\{\begin{array}{l@{\qquad}l}\displaystyle
    P(\partial_x,\partial_t)u=0,&x\in\mathbb{R}^n,\quad{t}\in\mathbb{R}^1.\\\displaystyle
    \left.\frac{\partial^ru}{\partial{t^r}}\right|_{t=0}=\varphi_r(x),&r=0,1,2,\ldots,2m-1.
    \end{array}\right.
\end{equation}
\begin{eqnarray}\label{5'}
  u(x,t) &=&
   \sum^m_{k=1}b_{2k}P(\partial_x)^{m-k}\sum^{2k-1}_{r=0}\frac{\partial^{2k-1-r}}{\partial t^{2k-1-r}}\int^t_0
   \frac{(t-\tau)^{2m-3}}{(2m-3)!}\nonumber\\
   & & \times\,\sum^m_{j=1}\frac{a_j^{2m-2}}{\prod^m_{i=1 \atop i\neq{j}}(a_j^2-a_i^2)}
   \frac{\sinh(\tau{a_j}P(\partial_x)^{1/2})}{a_jP(\partial_x)^{1/2}}\,\varphi_r(x)\,d\tau.
\end{eqnarray}

\textbf{Proof.} Considering initial conditions, the Laplace
transform of the Eq (\ref{5}) with respect to $t$ is
$$\sum^m_{k=0}b_{2k}P(\partial_x)^{m-k}(s^{2k}U(x,s)-\sum^{2k-1}_{r=0}s^{2k-1-r}\varphi_r(x))=0,$$
$$\prod^m_{i=1}(s^2-a_i^2P(\partial_x))U(x,s)-\sum^m_{k=1}b_{2k}P(\partial_x)^{m-k}\sum^{2k-1}_{r=0}s^{2k-1-r}\varphi_r(x)=0.$$
Where $U(x,s)=Lu(x,t)$. Let
$G_m(\partial_x,t)=L^{-1}[1/\prod^m_{i=1}(s^2-a_i^2P(\partial_x))]$,
by solving $U(x,s)$, we have its inverse Laplace transform:
\begin{eqnarray}\label{5''}
   u(x,t)&=& L^{-1}U(x,s)= \sum^m_{k=1}b_{2k}P(\partial_x)^{m-k}\sum^{2k-1}_{r=0}L^{-1}\frac{s^{2k-1-r}}{\prod^m_{i=1}(s^2-a_i^2P(\partial_x))}\,
   \varphi_r(x)\nonumber\\
   &=& \sum^m_{k=1}b_{2k}P(\partial_x)^{m-k}\sum^{2k-1}_{r=0}\frac{\partial^{2k-1-r}}{\partial{t^{2k-1-r}}}\,G_m(\partial_x,t)\varphi_r(x).
\end{eqnarray}

Now let us solve the $G_m(\partial_x,t)$. Considering initial
conditions, the Laplace transform of the Eq (\ref{I2}) with respect
to $t$ is
$$\prod^m_{i=1}(s^2-a_i^2P(\partial_x))U(x,s)=F(x,s),\quad{F(x,s)}=Lf(x,t).$$
By solving $U(x,s)$ and using the convolution theorem, we have its
inverse Laplace transform:
$$u(x,t)=L^{-1}U(x,s)=L^{-1}\frac{1}{\prod^m_{i=1}(s^2-a_i^2P(\partial_x))}F(x,s)=G_m(\partial_x,t)*f(x,t).$$
By comparing (\ref{I2'}) with $u(x,t)=G_m(\partial_x,t)*f(x,t)$, we
have the expression of the abstract operator $G_m(\partial_x,t)$:
\begin{equation}\label{6}
G_m(\partial_x,t)=\int^t_0\frac{(t-\tau)^{2m-3}}{(2m-3)!}
\sum^m_{j=1}\frac{a_j^{2m-2}}{\prod^m_{i=1 \atop
i\neq{j}}(a_j^2-a_i^2)}\frac{\sinh(\tau{a_j}P(\partial_x)^{1/2})}{a_jP(\partial_x)^{1/2}}\,d\tau.
\end{equation}
Applying (\ref{6}) to (\ref{5''}), thus the Corollary 2 is proved.

\section{Examples}

\textbf{Theorem 3.} Let $a_1,a_2,\ldots,a_m$ be arbitrary real roots
different from each other, which satisfy
$\sum^m_{k=0}b_{2k}x^{2k}=\prod^m_{i=1}(x^2-a_i^2)$, and
$P(\partial_x,\partial_t)$ be a partial differential operator
defined by
$$P(\partial_x,\partial_t)=\sum^m_{k=0}b_{2k}\Delta_n^{m-k}\frac{\partial^{2k}}{\partial{t^{2k}}},
\quad{x}\in\mathbb{R}^n,\;t\in\mathbb{R}^1,\;1<m\in\mathbb{N}.$$
Where $\Delta_n$ is an n-dimensional Laplacian, and
$n-2=2\nu+1,\;\nu\in\mathbb{N}_0$, then
$\forall{f(x,t)}\in{C^\infty}(\mathbb{R}_x^n\times\mathbb{R}_t^1),\;\varphi_r(x)\in{C^\infty}(\mathbb{R}^n)$,
we have
\begin{equation}\label{7}
    \left\{\begin{array}{l@{\qquad}l}\displaystyle
    P(\partial_x,\partial_t)u=f(x,t),&x\in\mathbb{R}^n,\quad{t}\in\mathbb{R}^1.\\\displaystyle
    \left.\frac{\partial^ru}{\partial{t^r}}\right|_{t=0}=\varphi_r(x),&r=0,1,2,\ldots,2m-1.
    \end{array}\right.
\end{equation}
\begin{eqnarray}\label{7'}
  u(x,t) &=&
  \int^t_0d\tau\int^{t-\tau}_0d\tau'\frac{(t-\tau-\tau')^{2m-3}}{(2m-3)!}\sum^m_{j=1}\frac{a_j^{2m-2}}{\prod^m_{i=1 \atop i\neq{j}}(a_j^2-a_i^2)}\nonumber\\
  & &
  \times\,[\tau'\underbrace{\int^{\tau'}_0\tau'd\tau'\cdots}_\nu\int^{\tau'}_0\frac{(a_j^{2}\Delta_n)^{\nu}}{S'_{n,j}}\int_{S'_{n,j}}f(\xi',\tau)\,dS'_{n,j}\,\tau'd\tau'
  +\sum^{\nu-1}_{l=0}\frac{a_j^{2l}\tau'^{2l+1}}{(2l+1)!}\Delta_n^lf(x,\tau)]\nonumber\\
   & & +\,\sum^m_{k=1}b_{2k}P(\partial_x)^{m-k}\sum^{2k-1}_{r=0}\frac{\partial^{2k-1-r}}{\partial t^{2k-1-r}}\int^t_0d\tau
   \frac{(t-\tau)^{2m-3}}{(2m-3)!}\sum^m_{j=1}\frac{a_j^{2m-2}}{\prod^m_{i=1 \atop i\neq{j}}(a_j^2-a_i^2)}\nonumber\\
   & &
   \times\,[\tau\underbrace{\int^{\tau}_0\tau{d\tau}\cdots}_\nu\int^{\tau}_0\frac{(a_j^{2}\Delta_n)^{\nu}}{S_{n,j}}\int_{S_{n,j}}\varphi_r(\xi)\,dS_{n,j}\,\tau{d\tau}
   +\sum^{\nu-1}_{l=0}\frac{a_j^{2l}\tau^{2l+1}}{(2l+1)!}\Delta_n^l\varphi_r(x)].
\end{eqnarray}
Where $S'_{n,j}=2(2\pi)^{\nu+1}(a_j\tau')^{n-1}$,
$S_{n,j}=2(2\pi)^{\nu+1}(a_j\tau)^{n-1}$, and $\xi'\in\mathbb{R}_n$
is the integral variable. The integral is on the hypersphere
$(\xi'_1-x_1)^2+(\xi'_2-x_2)^2+\cdots+(\xi'_n-x_n)^2=(a_j\tau')^2$,
and $dS'_{n,j}$ is its surface element. $\xi\in\mathbb{R}_n$ is the
integral variable on the hypersphere
$(\xi_1-x_1)^2+(\xi_2-x_2)^2+\cdots+(\xi_n-x_n)^2=(a_j\tau)^2$, and
$dS_{n,j}$ is its surface element.

\textbf{Proof.} According to (30) in reference \cite{bi10}, we have:
\begin{eqnarray}\label{8}
\frac{\sinh\left(ta_j\Delta_n{}^{1/2}\right)}{a_j\Delta_n{}^{1/2}}f(x)&=&
t\underbrace{\int^t_0tdt\cdots}_\nu\int^t_0\frac{(a_j^{2}\Delta_n)^{\nu}}{S_{n,j}}\int_{S_{n,j}}f(\xi)\,dS_{n,j}\,tdt\nonumber\\
   & &
   +\,\sum^{\nu-1}_{l=0}\frac{a_j^{2l}t^{2l+1}}{(2l+1)!}\Delta_n^lf(x),\quad\nu=\frac{n-3}{2}.
\end{eqnarray}
Where $\Delta_n$ is an n-dimensional Laplacian, $n-2=2\nu+1$,
$S_{n,j}=2(2\pi)^{\nu+1}(a_jt)^{n-1}$. $\xi\in\mathbb{R}_n$ is the
integral variable on the hypersphere
$(\xi_1-x_1)^2+(\xi_2-x_2)^2+\cdots+(\xi_n-x_n)^2=(a_jt)^2$, and
$dS_{n,j}$ is its surface element.

In Theorem 2, let $P(\partial_x)=\Delta_n$, then Theorem 3 is proved
by the substitution of (\ref{8}).

Similarly, we can easily obtain explicit solutions of the Cauchy
problem of more complex partial differential equations. For the
initial-boundary value problem, the operator $P(\partial_x)$ must
have the characteristic function related to boundary conditions, in
order to expand the known function $f(x,\tau),\;\varphi_r(x)$ in
(\ref{1'}) or (\ref{2'}) by using the characteristic function of
$P(\partial_x)$. Therefore, if the determination of the
characteristic function can be ascribed to the Sturm-Liouville
problem of given boundary conditions, then this initial-boundary
value problem is solvable.

$P(\partial_x)$ in Theorem BI4, Theorem 1 and Theorem 2 can be
variable-coefficient partial differential operators. For instance,
if $P(\partial_x)$ is a self-adjoint operator defined in a Hilbert
space, then the abstract operator
$$\frac{\sinh(tP(\partial_x)^{1/2})}{P(\partial_x)^{1/2}}\quad\mbox{and}\quad\cosh(tP(\partial_x)^{1/2})=\frac{\partial}{\partial{t}}\,
\frac{\sinh(tP(\partial_x)^{1/2})}{P(\partial_x)^{1/2}}$$ can act on
the Hilbert space, which also is a bounded operator in the Hilbert
space. In this case, we can attach proper boundary conditions to the
initial value problems in (\ref{I4}), (\ref{1}) and (\ref{2}).
Therefore, the given function $f(x,t),\varphi_r(x)$ becomes a
function with boundary conditions, and can be expressed in a Hilbert
space within the given domain. In order to solve the corresponding
initial-boundary value problem, we need to solve the characteristic
value problem of $P(\partial_x)$ under given boundary conditions to
determine a set of orthogonal functions, which generates a linear
manifold of a Hilbert space, thus $f(x,t),\varphi_r(x)$ can be
expressed in the Hilbert space.

\end{CJK*}
\end{document}